\documentclass[11pt]{article}
\usepackage{amsfonts}
\usepackage{amsfonts}
\usepackage{indentfirst,latexsym,bm,amsmath,amsthm}

\newtheorem{theorem}{Theorem}[section]
\newtheorem{prop}{Proposition}[section]
\newtheorem{definition}{Defintion}[section]
\renewcommand{\theequation}{\arabic{section}.\arabic{equation}}
\begin{document}
\title{\bf Existence and uniqueness of invariant measures for SPDEs with two reflecting walls
\footnote{$^1$ School of Mathematics, University of Manchester, Oxford Road, Manchester M13 9PL,
England, U.K. Email: tzhang@maths.man.ac.uk}}
\author{Juan Yang$^1$, Tusheng Zhang$^1$}

\date{}
\maketitle

\begin{abstract}
In this article, we study stochastic partial differential equations with two reflecting walls
$h^1$ and $h^2$, driven by space-time white noise with non-constant diffusion coefficients under periodic boundary conditions. The existence and uniqueness of invariant measures
is established under appropriate conditions. The strong Feller property is also obtained.\\

\noindent{\it Key Words:} stochastic partial differential equations with two reflecting walls; white noise; heat equation; invariant measures; coupling; strong Feller property.\\

\noindent{\it MSC:} Primary 60H15; Secondary 60J35
\end{abstract}

\section{Introduction} \setcounter{equation}{0}
\renewcommand{\theequation}{1.\arabic{equation}}

Consider the following stochastic partial differential equations (SPDEs) with two reflecting walls
\begin{equation}\label{1.1}
\left \{\begin{array}{ll}
        \frac{\partial{u(x,t)}}{\partial{t}}=\frac{\partial^2{u(x,t)}}{\partial{x^2}}+f\big(u(x,t)\big)+\sigma\big(u(x,t)\big)\dot{W}(x,t)\\         \ \ \ \ \ \ \ \ \ \ \ \   +\eta(x,t)-\xi(x,t);\\

        u(x,0)=u_0(x)\in{C}(S^1);\\
        h^1(x)\leq{u}(x,t)\leq{h^2}(x),\ \ {\rm for}\ (x,t)\in{Q}.
\end{array}
\right.
\end{equation}\\
${Q}:=S^1\times\mathbb{R}_+$, $S^1:=\mathbb{R}(mod 2\pi)$, or $\{e^{i\theta};\ \theta \in \mathbb{R}\}$ denotes a circular ring and the random field $W(x,t):=W(\{e^{i\theta};\ 0\leq \theta \leq x\}\times [0,t] )$ is a regular Brownian sheet defined on a filtered probability space
$(\Omega,P,\mathcal {F}; \mathcal {F}_t)$.
The random measures $\xi$ and $\eta$ are added to equation (\ref{1.1}) to
prevent the solution from leaving the interval $[h^1,\ h^2]$.

We assume that the reflecting walls $h^1(x),\ h^2(x)$ are continuous
functions satisfying
\\
(H1) $h^1(x)<h^2(x)$ for $x\in S^1$; \\
(H2)
$\frac{\partial^2{h^i}}{\partial{x^2}}\in
L^2(S^1)$, where $\frac{\partial^2{}}{\partial{x^2}}$ is interpreted in a distributional sense.\\
We also assume that the coefficients: $f, \sigma:
\mathbb{R}\rightarrow\mathbb{R}$
satisfy\\
(F1) there exists $L>0$ such that
$$|{f}(z_1)-f(z_2)|+|\sigma(z_1)-\sigma(z_2)|\leq{L}|z_1-z_2|,\ z_1,\ z_2\in \mathbb{R};$$
The following is the definition of a solution of a SPDE with two reflecting walls $h^1,\ h^2$.
\begin{definition} A triplet
$(u,\eta,\xi)$  is a solution to the
SPDE
(\ref{1.1}) if\\
(i) $u=\{u(x,t);(x,t)\in{Q}\}$ is a continuous, adapted random field
(i.e., $u(x,t)$ is $\mathcal {F}_t$-measurable $\forall$ $t\geq0,
x\in S^1$) satisfying $h^1(x)\leq{u}(x,t)\leq{h^2}(x)$,
a.s;\\
(ii) $\eta(dx,dt)$ and $\xi(dx,dt)$ are positive and adapted (i.e.
$\eta(B)$ and $\xi(B)$ is $\mathcal {F}_t$-measurable if
$B\subset S^1\times[0,t]$) random measures on $Q$
satisfying
$$\eta\big(S^1\times[0,T]\big)<\infty,\ \xi\big(S^1\times[0,T]\big)<\infty$$
for $T>0$;\\
(iii) for all $t\geq0$ and $\phi\in{C^\infty}(S^1)$ we
have
\begin{eqnarray}
&&\big(u(t),\phi\big)-\int_0^t(u(s),\phi^{''})ds-\int_0^t\big(f(u(s)),\phi\big)ds-
\int_0^t\int_{S^1}\phi(x)\sigma(u(x,s))W(dx,ds)\nonumber\\
&=&\big(u_0,\phi(x)\big)+\int_0^t\int_{S^1}\phi(x)\eta(dx,ds)-\int_0^t\int_{S^1}\phi(x)\xi(dx,ds),\ a.s,
\end{eqnarray}
where $(,)$ denotes the inner product in $L^2({S^1})$ and $u(t)$ denotes $u(\cdot,t)$;\\
(iv)$$\int_Q\big(u(x,t)-h^1(x)\big)\eta(dx,dt)=\int_Q\big(h^2(x)-u(x,t)\big)\xi(dx,dt)=0.$$
\end{definition}
\vskip 0.3cm
The existence and uniqueness of the solution of equation (\ref{1.1}) is established in \cite{ZY}, see also \cite{XZ} for SPDEs with one reflecting barrier. SPDEs with reflection were first studied by Nualart and Pardoux in \cite{NP}. Interesting properties were obtained in \cite{ZL}.

The aim of this paper is to establish the existence and uniqueness of invariant measures, as well as
the strong Feller property of fully non-linear SPDEs with two reflecting walls (\ref{1.1}).

For SPDEs without reflection, the existence and uniqueness of invariant measures has been studied by many people, see Sowers \cite{SL}, Mueller \cite{MC}, Peszat and Zabczyk \cite{PZ}, Da Prato and Zabczyk \cite{DZ}. For SPDEs with reflection, when the diffusion coefficient $\sigma$ is a constant, existence and uniqueness of invariant measures was obtained by  Otobe \cite{OI}, \cite{OS}. The strong Feller property of SPDEs has been studied by several authors, see
Peszat and Zabczyk \cite{PZ}, Da Prato and Zabczyk \cite{DZ}. The strong Feller property of SPDEs with reflection at $0$ was first proved in \cite{ZW}.

For the existence of invariant measures, our approach is to use Krylov-Bogolyubov theorem. To this end, the continuity of the solution with respect to the solutions of some random obstacle problems plays an important role. For the uniqueness, we adapted a coupling method used by Mueller \cite{MC}. Because of the reflection, we need to establish a kind of uniform coupling for approximating solutions. The strong Feller property of SPDEs with two reflecting walls will be obtained in a similar way as that  that  in Zhang \cite{ZW}.
\vskip 0.3cm

The rest of the paper is organized as follows. In Section 2, we give the proof of the existence and uniqueness of invariant measures. Section 3 establishes the strong Feller property.

\section{Existence and Uniqueness of Invariant Measures}
\setcounter{equation}{0}
\renewcommand{\theequation}{2.\arabic{equation}}

Denote by $\mathcal{B}({C}({S^1}))$ the $\sigma$-field
of all Borel subsets of ${C}({S^1})$ and by $\mathcal{M}({C}({S^1}))$ the set of all probability measures defined
on $({C}({S^1}),\ \mathcal{B}({C}({S^1})))$. We denote by $u(x,t,u_0)$ the solution of equation (\ref{1.1}) and by $P_t(u_0,\cdot)$ the corresponding transition function
$$P_t(u_0,\Gamma)=P(u(\cdot,t,u_0)\in\Gamma),\ \Gamma\in\mathcal{B}({C}({S^1})),\ t>0,$$
where $u_0$ is the initial condition. For
$\mu\in\mathcal{M}({C}({S^1}))$ we set $$P_t^*\mu(\Gamma)=\int_{C({S^1})}P_t(x,\Gamma)\mu(dx),$$ where $t\geq0,\ \Gamma\in\mathcal{B}({C}({S^1})).$
\begin{definition} A probability measure $\mu\in\mathcal{M}({C}({S^1}))$ is said to be invariant or stationary
with respect to $P_t$, $t\geq0$, if and only if $P_t^*\mu=\mu$ for each $t\geq0$.
\end{definition}

The initial condition $u_0(x)$ satisfies\\
(F2) $u_0(x)\in{}C({S^1})$ satisfy
$h^1(x)\leq{u_0}(x)\leq{h^2}(x),$ for $x\in S^1$.

\begin{theorem} Suppose the hypotheses (H1)-(H2), (F1)-(F2) hold.
Then there exists an invariant measure to equation (\ref{1.1}) on $C({S^1})$.
\end{theorem}
\noindent\textbf{Proof.}\quad
According to Krylov-Bogolyubov theorem (see \cite{DZ}), if the family $\{ P_t(u_0, \cdot);\ t\geq1\}$
is tight, then there exists an invariant measure for equation (1.1).  So we need to show that for any $\varepsilon>0$
there is a compact set $K\subset{C}({S^1})$ such that
\begin{eqnarray*}
{P}(u(t)\in{K})\geq1-\varepsilon,\ \  {\rm for\ any}\ t\geq1.
\end{eqnarray*}
where $u(t)=u(t,u_0)=u(\cdot,t,u_0)$.
On the other hand, for any $t\geq1$, we have by the Markov property
\begin{eqnarray}
{P}(u(t)\in{K})&=&\mathbb{E}\big(P_1(u(t-1), K)\big).
\end{eqnarray}
Thus it is enough to show ${P}\big(u(1,u(t-1))\in{K}\big)\geq1-\varepsilon$, for any $t\geq1$.
As $h^1(\cdot)\leq{}u(t-1)(\cdot)\leq{}h^2(\cdot)$, it suffices to find a compact subset $K\subset{C}({S^1})$ such that
\begin{eqnarray}\label{0.01}
{P}_1(g, K)\geq1-\varepsilon, \ \ for\ all\ g\in{C}({S^1})\ with\ h^1\leq{}g\leq{}h^2.
\end{eqnarray}
Put
\begin{eqnarray}\label{2.2}
        v(x,t,g)&=&\int_0^t\int_{S^1}G_{t-s}(x,y)f(u(y,s,g))dyds \nonumber\\
        &&+\int_0^t\int_{S^1}G_{t-s}(x,y)\sigma(u(y,s,g))W(dy,ds),
\end{eqnarray}
where $G_t(x,y)$ is the Green's function of the heat equation on $S^1$. Then $u$ can be written in the form(see \cite{NP}, \cite{DP1} and \cite{WA})
\begin{eqnarray*}
        u(x,t,g)-\int_{S^1}G_t(x,y)g(y)dy&=&v(x,t,g)+\int_0^t\int_{S^1}G_{t-s}(x,y)\eta(g)(dx,dt)\\
        &&-\int_0^t\int_{S^1}G_{t-s}(x,y)\xi(g)(dx,dt),
\end{eqnarray*}
where $\eta(g)$, $\xi(g)$ indicates the dependence of the random measures on the initial condition $g$.
Put
$$\bar{u}(x,t,g)=u(x,t,g)-\int_{S^1}G_t(x,y)g(y)dy$$
Then $(\bar{u}, \eta, \xi )$ solves a random obstacle problem.
From the relationship between $\bar{u}$ and $v$ proved in Theorem 4.1 in \cite{ZY}, we have the following inequality
$$\|\bar{u}(g)-\bar{u}(\hat{g})\|_\infty^1\leq 2\|v(g)-v(\hat{g})\|_\infty^1,$$
where $\|\omega\|_\infty^1:=\underset{x\in{S^1},t\in[0,1]}{\sup}|\omega(x,t)|.$
So $\bar{u}$ is a continuous functional of $v$ and denoted by $u=\Phi(v)$, where
 $\Phi(\cdot):\ {C}({S^1}\times{[0,1]})\rightarrow{C}({S^1}\times{[0,1]})$
is continuous. In particular, $\bar{u}(\cdot,1,g)$ is also a continuous functional of $v$, from ${C}({S^1}\times[0,1])$ to ${C}({S^1})$. We denote this functional by $\Phi_1$, i.e. $\bar{u}(\cdot,1,g)=\Phi_1(v(\cdot,g))$, where $v(\cdot,g)=v(\cdot,\cdot,g)$.
If $K''$ is a compact subset of ${C}({S^1}\times[0,1]))$,
then $K'=\Phi_1(K'')$ is a compact subset in $C({S^1})$ and
\begin{eqnarray}\label{0.02}
        P(\bar{u}(\cdot,1,g)\in{}K')
        &=&P(\bar{u}(\cdot,1,g)\in{}\Phi_1(K''))\nonumber \\
        &\geq&P(v(\cdot,g)\in{}K'').
\end{eqnarray}
Next, we want to find a compact set $K''(\subset{C}({S^1}\times[0,1])$ such that
\begin{equation}\label{0.03}
{P}(v(\cdot, g)\in{K''})\geq1-\varepsilon,  \ \ for\ all\ g\in{C}({S^1})\ with\ h^1\leq{}g\leq{}h^2.
\end{equation}

For $0<\alpha<\frac{1}{4}$ and $\kappa>0$, from Proposition A.1 in \cite{SL1} and
using a similar proof to that of Corollary 3.4 in \cite{WA},
there exists a random variable $Y(g)$ such that with probability one, for all $x,y\in{S^1}$ and $s,t\in(0,1]$,
\begin{eqnarray}\label{2.6}
|v(x,t,g)-v(y,s,g)|\leq{Y}(g)(d((x,t),(y,s)))^{\alpha-\kappa}\ and\ \mathbb{E}(Y(g))^{\frac{1}{\kappa}}\leq{C_0}, \
\end{eqnarray}
where $d((x,t),(y,s)):=\big(r^2(x,y)+(t-s)^2\big)^{\frac{1}{2}}$ with $r(x,y)$ the length of the shortest arc of $S^1$ connecting $x$ with $y$ and $C_0$ is independent of $g$.

Define
\begin{eqnarray*}
    \|v\|_\alpha=\sup&\{&\frac{|v(x,t)-v(y,s)|}{d^\alpha((x,t),(y,s))};\\
      &&(x,t),(y,s)\in{S^1}\times[0,1], (x,t)\neq(y,s)\}, \  for \ \alpha<\frac{1}{4}.
\end{eqnarray*}

By the Arzela-Ascoli theorem, for all $r>0$, $K_r:=\{v; \|v\|_\alpha\leq{r}\}$ is a compact subset of
$C({S^1}\times[0,1])$. In view of (\ref{2.6}), we see that for given $\varepsilon>0$, there exists $r_0$ such that
$${P}(v(\cdot,g)\in K_{r_0}^c)\leq\varepsilon,  \ \ for\ all\ g\ with\ h^1\leq{}g\leq{}h^2.$$
Choosing $K^{\prime\prime}=K_{r_0}$, we obtain (\ref{0.03}). Hence
$P(\bar{u}(\cdot,1,g)\in{}K')
        \geq 1-\varepsilon$   for all $g\in{C}({S^1})$  with $ h^1\leq{}g\leq{}h^2$. On the other hand, it is easy to see that there is a compact subset $K_0\subset C(S^1)$ such that
$$ \{ \int_{S^1}G_1(x,y)g(y)dy; \quad h^1\leq{}g\leq{}h^2 \}\subset K_0$$
Define $K=K^{\prime}+K_0$. We have
$$P_1(g,K)=P({u}(\cdot,1,g)\in{}K )\geq P(\bar{u}(\cdot,1,g)\in{}K')\geq 1-\varepsilon,$$
for all $g\in{C}({S^1})$  with $ h^1\leq{}g\leq{}h^2$. This finishes the proof.
\hfill$\Box$ \\

For the uniqueness of invariant measures, we need the following proposition. For simplicity, we put $u(x,t)=u(x,t,u_0)$.
\begin{prop} Under the assumption in Theorem 2.1, for any $p\geq1$, $T>0$,
$\underset{\varepsilon,\delta}{\sup}\mathbb{E}(\|u^{\varepsilon,\delta}\|_\infty^T)^p<\infty$ and
$u^{\varepsilon,\delta}$ converges uniformly on ${S^1}\times[0,T]$ to $u$ as $\varepsilon,\delta\rightarrow0$ a.s, where $u$, $u^{\varepsilon,\delta}$ are the solutions of equation (\ref{1.1}) and the penalized SPDEs
\begin{eqnarray*}
\left \{\begin{array}{ll}
       \frac{\partial{u^{\varepsilon,\delta}(x,t)}}{\partial{t}}=\frac{\partial^2{u^{\varepsilon,\delta}(x,t)}}{\partial{x^2}}+f(u^{\varepsilon,\delta}(x,t))
       +\sigma(u^{\varepsilon,\delta}(x,t))\dot{W}(x,t)\\
        \ \ \ \ \ \ \ \ \ \ \ \ \ \ +\frac{1}{\delta}(u^{\varepsilon,\delta}(x,t)-h^1(x))^- -\frac{1}{\varepsilon}(u^{\varepsilon,\delta}(x,t)-h^2(x))^+; \\
       u^{\varepsilon,\delta}(x,0)=u_0(x).
\end{array}
\right.
\end{eqnarray*}
\end{prop}
\noindent\textbf{Proof.}\quad Let $v^{\varepsilon,\delta}$ be the solution of equation
\begin{eqnarray}\label{2.7}
\left \{\begin{array}{ll}
       \frac{\partial{v^{\varepsilon,\delta}(x,t)}}{\partial{t}}=\frac{\partial^2{v^{\varepsilon,\delta}(x,t)}}{\partial{x^2}}+f(u^{\varepsilon,\delta}(x,t))
       +\sigma(u^{\varepsilon,\delta}(x,t))\dot{W}(x,t);\\
       v^{\varepsilon,\delta}(x,0)=u_0(x).
\end{array}
\right.
\end{eqnarray}
Set $\bar{\Phi}^{\varepsilon,\delta}(t)=\underset{s\leq{t},y\in{S^1}}{\sup}(v^{\varepsilon,\delta}(y,s)-h^2(y))^+$.
Note that $\bar{\Phi}^{\varepsilon,\delta}(t)$ is increasing w.r.t. $t$ and $v^{\varepsilon,\delta}-\bar{\Phi}^{\varepsilon,\delta}\leq{h^2}$.
$\bar{z}^{\varepsilon,\delta}(x,t):=v^{\varepsilon,\delta}(x,t)-\bar{\Phi}^{\varepsilon,\delta}(t)-u^{\varepsilon,\delta}(x,t)$ is a solution of equation
\begin{eqnarray}\label{2.8}
\left \{\begin{array}{ll}
       \frac{\partial{\bar{z}^{\varepsilon,\delta}}}{\partial{t}}+\frac{\partial{\bar{\Phi}^{\varepsilon,\delta}}}{\partial{t}}=\frac{\partial^2{\bar{z}^{\varepsilon,\delta}}}{\partial{x^2}}
       -\frac{1}{\delta}(u^{\varepsilon,\delta}-h^1)^- +\frac{1}{\varepsilon}(u^{\varepsilon,\delta}-h^2)^+; \\
       \bar{z}^{\varepsilon,\delta}(x,0)=0.
\end{array}
\right.
\end{eqnarray}
Multiplying (\ref{2.8}) by $(\bar{z}^{\varepsilon,\delta})^+$ and using $((u^{\varepsilon,\delta}-h^2)^+,\ (\bar{z}^{\varepsilon,\delta})^+)=0$ we get $(\bar{z}^{\varepsilon,\delta})^+=0$. Hence, $$u^{\varepsilon,\delta}\geq{v}^{\varepsilon,\delta}-\bar{\Phi}^{\varepsilon,\delta}.$$

Similarly, setting $\bar{z}^{\varepsilon,\delta}(x,t)=u^{\varepsilon,\delta}(x,t)-v^{\varepsilon,\delta}(x,t)-\underset{s\leq{t},y\in{S^1}}{\sup}(v^{\varepsilon,\delta}(y,s)-h^1(s))^-,$
we can show that
$$u^{\varepsilon,\delta}\leq{v}^{\varepsilon,\delta}+\underset{s\leq{t},y\in{S^1}}{\sup}(v^{\varepsilon,\delta}-h^1)^-.$$

As $\underset{\varepsilon,\delta}{\sup}\mathbb{E}(\|v^{\varepsilon,\delta}\|_\infty^T)^p<\infty$, the above two inequalities implies $$\underset{\varepsilon,\delta}{\sup}\mathbb{E}(\|u^{\varepsilon,\delta}\|_\infty^T)^p<\infty.$$
Since $u^{\varepsilon,\delta}$ is increasing in $\delta$ by the comparison theorem of SPDEs (see \cite{DP1}), we can show $u^\varepsilon:=\underset{\delta\downarrow0}{\lim}u^{\varepsilon,\delta}$ exists a.s. and $u^\varepsilon$ solves
\begin{eqnarray}\label{2.9}
\left \{\begin{array}{ll}
       \frac{\partial{u^{\varepsilon}(x,t)}}{\partial{t}}=\frac{\partial^2{u^{\varepsilon}(x,t)}}{\partial{x^2}}+f(u^{\varepsilon}(x,t))
       +\sigma(u^{\varepsilon}(x,t))\dot{W}(x,t)\\ \ \ \ \ \ \ \ \ \ \ \ \ +\eta^{\varepsilon}(x,t) -\frac{1}{\varepsilon}(u^{\varepsilon}(x,t)-h^2(x))^+; \\
       u^{\varepsilon}(x,0)\geq{h}^1(x);\\
       u^{\varepsilon}(x,0)=u_0(x),
\end{array}
\right.
\end{eqnarray}
where $\eta^\varepsilon(dx,dt):=\underset{\delta\downarrow0}{\lim}\frac{(u^{\varepsilon,\delta}(x,t)-h^1(x))^+}{\delta}dxdt$. Also, by comparison, we know that $u^{{\varepsilon}}$ is decreasing as $\varepsilon\downarrow0$.
Let $v^{\varepsilon}$ be the solution of equation (\ref{2.7}) replacing $u^{\varepsilon,\delta}$ by $u^{\varepsilon}$. Setting
$\bar{z}^{\varepsilon}(x,t)=u^{\varepsilon}(x,t)-v^{\varepsilon}(x,t)-\underset{s\leq{t},y\in{S^1}}{\sup}(v^{\varepsilon}(y,s)-h^1(y))^-,$
we can show $$u^{\varepsilon}\leq{v}^{\varepsilon}+\underset{s\leq{t},y\in{S^1}}{\sup}(v^{\varepsilon}-h^1)^-.$$ In addition, by the definition of $u^{{\varepsilon}}$, $u^{\varepsilon}\geq{h}^1$. Hence, $u:=\underset{\varepsilon\downarrow0}{\lim}u^\varepsilon=\underset{\varepsilon\downarrow0}{\lim}\underset{\delta\downarrow0}{\lim}u^{\varepsilon,\delta}$ exists a.s.

The continuity of $u$ can be proved similarly as in Theorem 4.1 in \cite{DP1}. The uniform convergence of $u^{\varepsilon,\delta}$ w.r.t. $(x,t)$ follows from Dini's theorem.

\hfill$\Box$\\

The following result is the uniqueness of invariant measures.
\begin{theorem} Under the assumptions in Theorem 2.1 and that $\sigma\geq{L_0}$ for some constant $L_0>0$,
there is a unique invariant measure for the equation (\ref{1.1}).
\end{theorem}
\noindent\textbf{Proof.}\quad We will adopt  the coupling method used in Mueller \cite{MC} to SPDEs with reflection.
Let $u^1(x,0)$ and $u^2(x,0)$ be two initial values having distributions given by two invariant probabilities $\mu_1$ and $\mu_2$. Then $ u^1(x,t)$ and $u^2(x,t)$ also have these distributions for any $t>0$.  Thus
$$Var(\mu_1-\mu_2)\leq P\big (\sup_{x\in {S^1}}|u^1(x,t)-u^2(x,t)|\not =0\big ).$$
Thus, for given two initial functions $u^1(x,0)$ and $u^2(x,0)$, it is sufficient to construct two coupled processes $u^1(x,t)$, $u^2(x,t)$ satisfying equation (\ref{1.1}), driven by different white noises on a probability space $({\Omega}, {\cal F}, {P})$, such that
  \begin{equation}\label{3.001}
  \underset{t\rightarrow \infty}{\lim}P\big (\sup_{x\in {S^1}}|u^1(x,t)-u^2(x,t)|\not =0\big )=0.
  \end{equation}

We first assume $u^1(x,0)\geq u^2(x,0)$, $x\in{S^1}$. We want to construct two independent space-time white noises ${W}_1(x,t),\ {W}_2(x,t)$ defined on a probability space $({\Omega}, {\cal F}, {P})$, and a solution  ${u}, {v}$  of the following SPDEs with two reflecting walls
\begin{eqnarray}\label{3.1}
\frac{\partial{u(x,t)}}{\partial{t}}&=&\frac{\partial^2{u(x,t)}}{\partial{x^2}}+f\big(u(x,t)\big)
+\sigma\big(u(x,t)\big)\dot{W}_1(x,t)\nonumber\\
                                     &&+\eta_1(x,t)-\xi_1(x,t),\nonumber\\
  \frac{\partial{v(x,t)}}{\partial{t}}&=&\frac{\partial^2{v(x,t)}}{\partial{x^2}}
  +f\big(v(x,t)\big)+\eta_2(x,t)-\xi_2(x,t)\nonumber\\
  &&+\sigma\big(v(x,t)\big)\big[(1-|u-v|\wedge 1)^{\frac{1}{2}}\dot{W}_1(x,t)+(|u-v|\wedge 1)^{\frac{1}{2}}\dot{W}_2(x,t)\big ],\nonumber\\
  u(x,0)&=&u^1(x,0), \quad v(x,0)=u^2(x,0).
  \end{eqnarray}
Note that the  coefficients in the second equation in (\ref{3.1}) is not Lipschitz. The existence of a solution of  equation (\ref{3.1}) is not automatic. In the following, using a similar method as that in the paper \cite{MC}, we will give a construction of a solution on some probability space. The construction will also be used to prove the successful coupling
$$\underset{t\rightarrow \infty}{\lim}P\big (\sup_{x\in {S^1}}|u(x,t)-v(x,t)|\not =0\big )=0.$$
For $0\leq z\leq 1$, set
  \begin{eqnarray*}
  f_n(z)&=&\big ( z+\frac{1}{n}\big )^{\frac{1}{2}}-\big ( \frac{1}{n}\big )^{\frac{1}{2}},\nonumber\\
  g_n(x)&=&\big (1-f_n(z)^2\big )^{\frac{1}{2}}.
  \end{eqnarray*}
  We have $f_n(z)^2+g_n(z)^2=1$ and that $f_n(z)\rightarrow z^{\frac{1}{2}}$, $g_n(z)\rightarrow (1-z)^{\frac{1}{2}}$ uniformly as $n\rightarrow \infty$, for $z\in {S^1}$.

  Let  $\dot{\overline{W}}_1(x,t), \ \dot{\overline{W}}_2(x,t)$ be two independent space-time white noises  defined on a probability space $(\overline{\Omega}, \overline{\cal F}, \overline{P})$. Let $\overline{u}, \overline{v}^n$ be the unique solution of the following SPDEs with two reflecting walls
  \begin{eqnarray}\label{3.2}
\frac{\partial{\overline{u}(x,t)}}{\partial{t}}&=&
\frac{\partial^2{\overline{u}(x,t)}}{\partial{x^2}}+f\big(\overline{u}(x,t)\big)
+\sigma\big(\overline{u}(x,t)\big)\dot{\overline{W}}_1(x,t)\nonumber\\
                                     &&+\overline{\eta}_1(x,t)-\overline{\xi}_1(x,t),\nonumber\\
  \frac{\partial{\overline{v}^n(x,t)}}{\partial{t}}&=&
  \frac{\partial^2{\overline{v}^n(x,t)}}{\partial{x^2}}+f\big(\overline{v}^n(x,t)\big)+\bar{\eta}^n_2(x,t)-\bar{\xi}^n_2(x,t)\nonumber\\
  &&+\sigma\big(\overline{v}^n(x,t)\big)\big[g_n(|\bar{u}-\overline{v}^n|\wedge 1)\dot{\overline{W}}_1(x,t)+f_n(|\overline{u}-\overline{v}^n|\wedge 1)\dot{\overline{W}}_2(x,t)\big ],\nonumber\\
  \overline{u}(x,0)&=&u^1(x,0), \quad \overline{v}(x,0)=u^2(x,0).
                                     \end{eqnarray}
  The existence and uniqueness of $(\overline{u}, \overline{v}^n)$ is guaranteed because of the Lipschitz continuity of the coefficients. Put
  \begin{eqnarray}\label{3.3}
        \hat{u}(x,t)&=&\int_{S^1}G_t(x,y)u^1(y,0)dy+\int_0^t\int_{S^1}G_{t-s}(x,y)
        f(\overline{u}(y,s))dyds \nonumber\\
        &&+\int_0^t\int_{S^1}G_{t-s}(x,y)\sigma(\overline{u}(y,s))\overline{W}_1(dy,ds)
\end{eqnarray}
and
\begin{eqnarray}\label{3.4}
        \hat{v}^n(x,t)&=&\int_{S^1}G_t(x,y)u^2(y,0)dy+\int_0^t\int_{S^1}G_{t-s}(x,y)
        f(\overline{v}^n(y,s))dyds \nonumber\\
        &&+\int_0^t\int_{S^1}G_{t-s}(x,y)\sigma(\overline{v}^n(y,s))W^n(dy,ds),
\end{eqnarray}
where
$$\dot{W}^n(x,t)=\big[g_n(|\overline{u}-\overline{v}^n|\wedge 1)\dot{\overline{W}}_1(x,t)+f_n(|\overline{u}-\overline{v}^n|\wedge 1)\dot{\overline{W}}_2(x,t)\big ]$$
is another space-time white noise on $(\overline{\Omega}, \overline{\cal F}, \overline{P})$.
 From the proof of Theorem 2.1, it is known that  there exists a continuous functional $\Phi$ from $C(S^1\times [0,T])$ into $C(S^1\times [0,T])$ (for ant $T>0$) such that $\bar{u}=\Phi(\hat{u})$ and $\bar{v}^n=\Phi(\hat{v}^n)$.
On the other hand, following the same proof of Lemma 3.1 in \cite{MC} it can be shown  that the sequence $\hat{u}, \hat{v}^n, n\geq 1$ is tight. As the images under  the continuous map $\Phi$, the vector $(\overline{u}, \overline{v}^n, \overline{W}_1, \overline{W}_2)$ is also tight. By Skorohod's representation theorem, there exist random fields $(u, v^n, W_1, W_2)$, $n\geq 1$ on some probability space $(\Omega, {\cal F}, P)$ such that $(u, v^n, W_1, W_2)$ has the same law as $(\overline{u}, \overline{v}^n, \overline{W}_1, \overline{W}_2)$ and that the following SPDEs with two reflecting walls hold
\begin{eqnarray}\label{3.5}
\frac{\partial{u(x,t)}}{\partial{t}}&=&\frac{\partial^2{{u}(x,t)}}{\partial{x^2}}+f\big(u(x,t)\big)
+\sigma\big({u}(x,t)\big)\dot{W}_1(x,t)\nonumber\\
                                     &&+\eta_1(x,t)-\xi_1(x,t),\nonumber\\
  \frac{\partial{{v}^n(x,t)}}{\partial{t}}&=&\frac{\partial^2{{v}^n(x,t)}}{\partial{x^2}}+f\big(v^n(x,t)\big)+\eta_2^n(x,t)-\xi^n_2(x,t)\nonumber\\
  &&+\sigma\big(v^n(x,t)\big)\big[g_n(|u-v^n|\wedge 1)\dot{W}_1(x,t)+f_n(|u-v^n|\wedge 1)\dot{W}_2(x,t)\big ],\nonumber\\
  u(x,0)&=&u^1(x,0), \quad v^n(x,0)=u^2(x,0).
\end{eqnarray}
Furthermore, $v^n\rightarrow v$ uniformly almost surely as $n\rightarrow \infty$. By a  similar proof as that of Theorem 4.1 in \cite{ZY} we can prove that
the limit $(u, v)$ satisfies the following SPDEs with two reflecting walls
\begin{eqnarray}\label{3.5}
\frac{\partial{u(x,t)}}{\partial{t}}&=&\frac{\partial^2{u(x,t)}}{\partial{x^2}}+f\big(u(x,t)\big)
+\sigma\big(u(x,t)\big)\dot{W}_1(x,t)\nonumber\\
                                     &&+\eta_1(x,t)-\xi_1(x,t),\nonumber\\
  \frac{\partial{v(x,t)}}{\partial{t}}&=&\frac{\partial^2{v(x,t)}}{\partial{x^2}}
  +f\big(v(x,t)\big)+\eta_2(x,t)-\xi_2(x,t)\nonumber\\
  &&+\sigma\big(v(x,t)\big)\big[(1-|u-v|\wedge 1)^{\frac{1}{2}}\dot{W}_1(x,t)+(|u-v|\wedge 1)^{\frac{1}{2}}\dot{W}_2(x,t)\big ],\nonumber\\
  u(x,0)&=&u^1(x,0), \quad v(x,0)=u^2(x,0).
\end{eqnarray}
The next step is to show that $u,v$ admits a successful coupling. To this end, consider the following approximating SPDEs
\begin{eqnarray}\label{3.0006}
\left \{\begin{array}{ll}
       \frac{\partial{{u}^{\varepsilon,\delta}}}{\partial{t}}=\frac{\partial^2{{u}^{\varepsilon,\delta}}}{\partial{x^2}}
       +f({u}^{\varepsilon,\delta})
       +\frac{1}{\delta}({u}^{\varepsilon,\delta}-h^1)^- -\frac{1}{\varepsilon}({u}^{\varepsilon,\delta}-h^2)^+
       +\sigma({u}^{\varepsilon,\delta})\dot{W}_1; \\
              \frac{\partial{{v}^{n, \varepsilon,\delta}}}{\partial{t}}=\frac{\partial^2{{v}^{n, \varepsilon,\delta}}}{\partial{x^2}}+f({v}^{n, \varepsilon,\delta})
       +\frac{1}{\delta}({v}^{n, \varepsilon,\delta}-h^1)^- -\frac{1}{\varepsilon}({v}^{n, \varepsilon,\delta}-h^2)^+ \\
        \ \ \ \ \ \ \ \ \ \ \ \
       +\sigma({v}^{n, \varepsilon,\delta})\big[g_n(|u^{\varepsilon,\delta}-v^{n, \varepsilon,\delta}|\wedge 1)\dot{W}_1(x,t)\\
       \ \ \ \ \ \ \ \ \ \ \ \
       +f_n(|u^{ \varepsilon,\delta}-v^{n, \varepsilon,\delta}|\wedge 1)\dot{W}_2(x,t)\big ];\\
              {u}^{\varepsilon,\delta}(x,0)={u}^1(x,0),\ {v}^{n, \varepsilon,\delta}(x,0)={u}^2(x,0).
\end{array}
\right.
\end{eqnarray}
We may and will assume that $f(u)$ is non-increasing. Otherwise, we consider $\tilde{u}:=e^{-Lt}u$, $\tilde{v}:=e^{-Lt}v$, where $L$ is the Lipschitz constant in (F1), which satisfy
\begin{eqnarray*}
\frac{\partial{\tilde{u}(x,t)}}{\partial{t}}&=&\frac{\partial^2{{\tilde{u}}(x,t)}}{\partial{x^2}}+e^{-Lt} f\big(e^{Lt}\tilde{u}(x,t)\big)- L\tilde{u}(x,t)\\
&& +e^{-Lt}\sigma\big(e^{Lt}{\tilde{u}}(x,t)\big)\dot{W}_1(x,t)
                                     +\eta_3(x,t)-\xi_3(x,t),\nonumber\\
  \frac{\partial{{\tilde{v}}^n(x,t)}}{\partial{t}}&=&\frac{\partial^2{{\tilde{v}}^n(x,t)}}{\partial{x^2}}+e^{-Lt} f\big(e^{Lt}\tilde{v}(x,t)\big)- L\tilde{v}(x,t) +\eta_4^n(x,t) -\xi^n_4(x,t)\\
&& +e^{-Lt}\sigma\big(e^{Lt}\tilde{v}^n(x,t)\big)\big[g_n(|e^{Lt}\tilde{u}-e^{Lt}\tilde{v}^n|\wedge 1)\dot{W}_1(x,t)\\&&+f_n(|e^{Lt}\tilde{u}-e^{Lt}\tilde{v}^n|\wedge 1)\dot{W}_2(x,t)\big ],\nonumber\\
  u(x,0)&=&u^1(x,0), \quad v^n(x,0)=u^2(x,0).
\end{eqnarray*}
The new drift $e^{-Lt}f(e^{Lt}x)-Lx$ is non-increasing. Also, if $\tilde{u}, \ \tilde{v}$ satisfy a successful coupling, so does $u, \ v$.
Note that all the coefficients in (\ref{3.0006}) are Lipschitz continuous. We can apply Proposition 2.1 to conclude that ${u}^{\varepsilon,\delta}(x,t)\rightarrow u(x,t)$, ${v}^{n, \varepsilon,\delta}(x,t)\rightarrow v^n(x,t)$ uniformly on ${S^1}\times [0,T]$ ( for any $T>0$) as $\varepsilon, \delta \rightarrow 0$. As $u^1(x,0)\geq u^2(x,0)$, as lemma 3.1 in \cite{MC}, we can show that ${u}^{\varepsilon,\delta}\geq{v}^{n, \varepsilon,\delta}$. Let
\begin{equation}\label{3.002}
U^{n,\varepsilon, \delta}(t)=\int_{S^1}({u}^{\varepsilon,\delta}(x,t)-{v}^{n,\varepsilon,\delta}(x,t))dx.
\end{equation}
It follows from the above equation that
\begin{equation}\label{3.003}
U^{n,\varepsilon, \delta}(t)=\int_{S^1}({u}_1(x,0)-{u}_2(x,0))dx+\int_0^t C^{n,\varepsilon, \delta}(s)ds+ M^{n,\varepsilon, \delta}(t),
\end{equation}
where
\begin{eqnarray*}
        C^{n,\varepsilon,\delta}(t)&=&\int_{S^1}\big \{f({u}^{\varepsilon,\delta})-f({v}^{n, \varepsilon,\delta})+\frac{1}{\delta}({u}^{\varepsilon,\delta}-h^1)^-(x,t)-\frac{1}{\delta}({v}^{n,\varepsilon,\delta}-h^1)^-(x,t)\\
&&-\big(\frac{1}{\varepsilon}({u}^{\varepsilon,\delta}-h^2)^+(x,t)-\frac{1}{\varepsilon}({v}^{n,\varepsilon,\delta}-h^2)^+(x,t)\big)\big\}dx\\
&\leq & 0,
\end{eqnarray*}

\begin{eqnarray*}
        M^{n,\varepsilon,\delta}(t)&=&\int_0^t\int_{S^1}\sigma({u}^{\varepsilon,\delta}(x,s)){W}_1(dx,ds) \\
&&-\int_0^t\int_{S^1}\sigma({v}^{n, \varepsilon,\delta}(x,s))g_n(|u^{\varepsilon,\delta}-v^{n, \varepsilon,\delta}|\wedge 1)\dot{W}_1(dx,ds)\\
&&- \int_0^t\int_{S^1}\sigma({v}^{n, \varepsilon,\delta}(x,s))f_n(|u^{ \varepsilon,\delta}-v^{n, \varepsilon,\delta}|\wedge 1)\dot{W}_2(dx,ds).
\end{eqnarray*}
Observe that
\begin{eqnarray}\label{3.004}
&&\lim_{\varepsilon, \delta\rightarrow 0}U^{n,\varepsilon, \delta}(t)\nonumber\\
&=&U^n(t):=\int_{S^1}({u}(x,t)-{v}^{n}(x,t))dx ,
\end{eqnarray}
and
\begin{eqnarray}\label{3.005}
        &&\lim_{\varepsilon, \delta\rightarrow 0}M^{n,\varepsilon,\delta}(t)\nonumber\\
        &=&M^n(t):=\int_0^t\int_{S^1}\sigma({u}(x,s)){W}_1(dx,ds) \nonumber \\
&&-\int_0^t\int_{S^1}\sigma({v}^{n}(x,s))g_n(|u-v^{n}|\wedge 1)\dot{W}_1(dx,ds)\nonumber\\
&&- \int_0^t\int_{S^1}\sigma({v}^{n}(x,s))f_n(|u-v^{n}|\wedge 1)\dot{W}_2(dx,ds).
\end{eqnarray}
Letting $\varepsilon, \delta \rightarrow 0$ in (\ref{3.003}) we see that
\begin{equation}\label{3.006}
U^{n}(t)=\int_{S^1}({u}_1(x,0)-{u}_2(x,0))dx+A^n(t)+ M^{n}(t),
\end{equation}
where $A^n(t)=\underset{\varepsilon,\delta\rightarrow 0}\lim\int_0^t C^{n,\varepsilon, \delta}(s)ds$ is a continuous, adapted  non-increasing process. Now, sending $n$ to $\infty$ we obtain
\begin{equation}\label{3.007}
U(t)=\int_{S^1}({u}_1(x,0)-{u}_2(x,0))dx+A(t)+ M(t),
\end{equation}
where
$$U(t)=\int_{S^1}({u}(x,t)-{v}(x,t))dx,$$
\begin{eqnarray*}
        M(t)&=&\int_0^t\int_{S^1}\sigma({u}(x,s)){W}_1(dx,ds) \\
&&-\int_0^t\int_{S^1}\sigma({v}(x,s))(1-|u-v|\wedge 1)^{\frac{1}{2}}\dot{W}_1(dx,ds)\\
&&- \int_0^t\int_{S^1}\sigma({v}(x,s))(|u-v|\wedge 1)^{\frac{1}{2}}\dot{W}_2(dx,ds),
\end{eqnarray*}
and $A(t)=\underset{n\rightarrow \infty}\lim A^n(t)$ a continuous, adapted  non-increasing process. The existence of the limits of $A^n$ follows from the existence of the limit of $U^n$ and $M^n$.
Now we can modify the proof in \cite{MC} to obtain the successful coupling of $u$ and $v$.
In view of the assumption on $\sigma$ and the boundedness of the walls $h^1,h^2$,  it is easy to verify that
\begin{equation}\label{3.008}
\frac{d<M>(t)}{dt}\geq C_0U(t)
\end{equation}
for some positive constant $C_0$. Thus, there exists a non-negative adapted process $V(t)$ such that
$$\frac{d<M>(t)}{dt}=U(t)V(t),\quad \quad V(t)\geq C_0. $$
Let
\begin{eqnarray}\label{3.009}
        \phi(t)&=&\int_0^tV(s)ds,\nonumber\\
        X(t)&=&U(\phi^{-1}(t)).
\end{eqnarray}
Then the time-changed process $X$ satisfies the following equation
\begin{equation}\label{3.0010}
X(t)=U(0)+\tilde{A}(t)+\int_0^tX^{\frac{1}{2}}(s)dB(s),
\end{equation}
where $B$ is a Brownian motion and $\tilde{A}$ is an adapted non-increasing process. Let $Y(t)=2X^{\frac{1}{2}}(t)$. Applying Ito's formula  (before $Y$ hits $0$) we obtain
\begin{equation}\label{3.0011}
Y(t)=Y(0)+2\int_0^t\frac{1}{Y(s)}d\tilde{A}(s)-\frac{1}{2}\int_0^t \frac{1}{Y(s)}ds+B(t).
\end{equation}
As $\tilde{A}$ is non-increasing, it follows that
\begin{equation}\label{3.0012}
0\leq Y(t)\leq Y(0)+B(t).
\end{equation}
The property of one dimensional Brownian motion implies that $Y$ hits $0$ with probability $1$. Hence

$$\underset{t\rightarrow \infty}{\lim}P\big (\sup_{x\in {S^1}}|u(x,t)-v(x,t)|\not =0\big )=0.$$

Next let us consider the general case, i.e. we do not assume $u^1(x,0)\geq u^2(x,0)$, $x\in{S^1}$. Consider a solution  $v$, $u^1$, $u^2$  of the following SPDEs with two reflecting walls
\begin{eqnarray*}
\frac{\partial{v(x,t)}}{\partial{t}}&=&\frac{\partial^2{v(x,t)}}{\partial{x^2}}+f\big(v(x,t)\big)
+\sigma\big(v(x,t)\big)\dot{W}_1(x,t)\nonumber\\
                                     &&+\eta_v(x,t)-\xi_v(x,t),\nonumber\\
  \frac{\partial{u^i(x,t)}}{\partial{t}}&=&\frac{\partial^2{u^i(x,t)}}{\partial{x^2}}
  +f\big(u^i(x,t)\big)+\eta_{u^i}(x,t)-\xi_{u^i}(x,t)\nonumber\\
  &&+\sigma\big(u^i(x,t)\big)\big[(1-|v-u^i|\wedge 1)^{\frac{1}{2}}\dot{W}_1(x,t)+(|v-u^i|\wedge 1)^{\frac{1}{2}}\dot{W}_2(x,t)\big ],\nonumber\\
 v(x,0)&=&\underset{i=1,2}{\max}\{u^i(x,0)\}.
 \end{eqnarray*}
By following the arguments in the first part, we have $$\underset{t\rightarrow \infty}{\lim}P\big (\sup_{x\in {S^1}}|v(x,t)-u^i(x,t)|\not =0\big )=0,\ i=1,2.$$
The inequality
$$0\leq \sup_{x\in {S^1}}|u^1(x,t)-u^2(x,t)|   \leq  \sum\limits_{i=1}^2  \big(\sup_{x\in {S^1}}|v(x,t)-u^i(x,t)|\big)$$
implies
$$\underset{t\rightarrow \infty}{\lim}P\big (\sup_{x\in {S^1}}|u^1(x,t)-u^2(x,t)|\not =0\big )=0.$$
\hfill$\Box$

\section{Strong Feller property}
\setcounter{equation}{0}
\renewcommand{\theequation}{3.\arabic{equation}}
In this section, we consider the strong Feller property of the solution of equation (\ref{1.1}). Let $H=L^2({S^1})$. If $\varphi\in{B}_b(H)$ (the Banach space of all real bounded Borel
functions, endowed with the sup norm), we define, for $x\in{S^1}$, $0\leq{}t\leq{}T$ and $g\in{}H$,
$$P_{t}\varphi(g)=\mathbb{E}\varphi(u(x,t,g)).$$

\begin{definition}
The family $\{P_{t}\}$ is called strong Feller if for arbitrary $\varphi\in{B}_b(H)$, the function $P_{t}\varphi(\cdot)$ is continuous for all $t>0$.
\end{definition}
\begin{theorem} Under the hypotheses (H1)-(H2), (F1)-(F3) and that $p_1\leq|\sigma(\cdot)|\leq{p_2}$ for some constants $p_1,\ p_2> 0$, then for any $T>0$ there exists a constant $C'_T$ such that for all
$\varphi\in{B}_b(H)$ and $t\in(0,T]$,
\begin{eqnarray}\label{3.01}
|P_{t}\varphi(u_0^1)-P_{t}\varphi(u_0^2)|&\leq&\frac{C'_T}{\sqrt{t}}\|\varphi\|_{\infty}|u_0^1-u_0^2|_H,
\end{eqnarray}
for $u_0^1,\ u_0^2\in H$ with
$h^1(x)\leq{u_0^1}(x),\ {u_0^2}(x)\leq{h^2}(x),$
where $\|\varphi\|_{\infty}=\underset{u_0}{\sup}|\varphi(u_0)|$. In particular, $P_t,\ t>0$, is strong Feller.
\end{theorem}
\noindent\textbf{Proof.}\quad
Choose a non-negative function $\phi\in{}C_0^\infty(R)$ with
$\int_R\phi(x)=1$ and denote $$f_n(\zeta)=n\int_R\phi\big(n(\zeta-y)\big)f(y)dy,$$
$$\sigma_n(\zeta)=n\int_R\phi\big(n(\zeta-y)\big)\sigma(y)dy,$$
$$k_n(\zeta,x)=n\int_R\phi\big(n(\zeta-y)\big)(y-h^1(x))^-dy,$$
$$l_n(\zeta,x)=n\int_R\phi\big(n(\zeta-y)\big)(y-h^2(x))^+dy.$$
So $f_n,\ \sigma_n,\ k_n,\ l_n$ are smooth w.r.t. $\zeta$.
Let
\begin{eqnarray*}
u_n^{\varepsilon,\delta}(x,t,u_0)&=&\int_{S^1}G_t(x,y)u_0(y)dy+\int_0^t\int_{S^1}G_{t-s}(x,y)f_n\big(u_n^{\varepsilon,\delta}(y,s,u_0)\big)dyds\\
&&+\int_0^t\int_{S^1}G_{t-s}(x,y)\sigma_n\big(u_n^{\varepsilon,\delta}(y,s,u_0)\big)W(dy,ds)\\
&&+\frac{1}{\delta}\int_0^t\int_{S^1}G_{t-s}(x,y)k_n\big(u_n^{\varepsilon,\delta}(y,s,u_0),y\big)dyds\\
&&-\frac{1}{\varepsilon}\int_0^t\int_{S^1}G_{t-s}(x,y)l_n\big(u_n^{\varepsilon,\delta}(y,s,u_0),y\big)dyds.
\end{eqnarray*}
Since $f_n(\zeta)\rightarrow{}f(\zeta)$, $\sigma_n(\zeta)\rightarrow{}\sigma(\zeta)$,
$k_n(\zeta,x)\rightarrow(\zeta-h^1(x))^-$ and $l_n(\zeta,x)\rightarrow(\zeta-h^2(x))^+$ as $n\rightarrow\infty,$
we can show that for any fixed $\varepsilon,\ \delta$ and $p\geq1$,
$$\underset{n\rightarrow\infty}{\lim}\underset{t\in[0,T]}{\sup}\mathbb{E}\big(|u_n^{\varepsilon,\delta}(t,\cdot,u_0)-u^{\varepsilon,\delta}(t,\cdot,u_0)|_H^p\big)=0.$$
By Lemma 7.1.5 in \cite{DZ} and Proposition 2.1, it is enough to prove that
there exists a constant $C'_T$, independent of $\varepsilon,\ \delta$ and $n$, such that
\begin{eqnarray}\label{3.02}
|P_{t}^{n,\varepsilon,\delta}\varphi(u_0^1)-P_{t}^{n,\varepsilon,\delta}\varphi(u_0^2)|&\leq&\frac{C'_T}{\sqrt{t}}\|\varphi\|_{\infty}|u_0^1-u_0^2|_H,
\end{eqnarray}
where $P_{t}^{n,\varepsilon,\delta}\varphi(u_0):=\mathbb{E}\big(\varphi(u_n^{\varepsilon,\delta}(\cdot,\cdot,u_0))\big)$ and $\ u_0^1$, $u_0^2\in{}H$.

From Theorem 5.4.1 in \cite{DZ}, $u_n^{\varepsilon,\delta}(\cdot, \cdot, u_0)$ is continuously differentiable w.r.t. $u_0$. Denote by
$X_n^{\varepsilon,\delta}(x,t):=\big(Du_n^{\varepsilon,\delta}(\cdot,\cdot,u_0)(\bar{u}_0)\big)(x,t)$ the directional derivative of $u_n^{\varepsilon,\delta}(\cdot,\cdot,u_0)$
at $u_0$ in the direction of $\bar{u}_0$ and it satisfies the mild form of a SPDE
\begin{eqnarray*}
X_n^{\varepsilon,\delta}(x,t)&=&\int_{S^1}G_t(x,y)\bar{u}_0(y)dy
+\int_0^t\int_{S^1}G_{t-s}(x,y)f'_n\big(u_n^{\varepsilon,\delta}(y,s,u_0)\big)X_n^{\varepsilon,\delta}(y,s)dyds\\
&&+\int_0^t\int_{S^1}G_{t-s}(x,y)\sigma'_n\big(u_n^{\varepsilon,\delta}(y,s,u_0)\big)X_n^{\varepsilon,\delta}(y,s)W(dy,ds)\\
&&+\frac{1}{\delta}\int_0^t\int_{S^1}G_{t-s}(x,y)\frac{\partial}{\partial{\zeta}}k_n\big(u_n^{\varepsilon,\delta}(y,s,u_0),y\big)X_n^{\varepsilon,\delta}(y,s)dyds\\
&&-\frac{1}{\varepsilon}\int_0^t\int_{S^1}G_{t-s}(x,y)\frac{\partial}{\partial{\zeta}}l_n\big(u_n^{\varepsilon,\delta}(y,s,u_0),y\big)X_n^{\varepsilon,\delta}(y,s)dyds.
\end{eqnarray*}
Since $\frac{\partial}{\partial{\zeta}}k_n\big(u_n^{\varepsilon,\delta}(y,s,u_0),y\big)\leq0$, $\frac{\partial}{\partial{\zeta}}l_n\big(u_n^{\varepsilon,\delta}(y,s,u_0),y\big)\geq0$,  we use the similar arguments as that in \cite{ZW} and to get
\begin{eqnarray*}
\underset{\varepsilon,\delta\geq0,t\in[0,T]}{\sup}\mathbb{E}\big(\int_{S^1}(X_n^{\varepsilon,\delta}(y,t))^2dy\big)&\leq&C|\bar{u}_0|_H^2,
\end{eqnarray*}
where $C$ is a constant.
By Elworthy-Li formula (Lemma 7.1.3 in \cite{DZ}), we obtain
\begin{eqnarray*}
\mid\langle{}DP_{t}\varphi(u_0),\bar{u}_0\rangle\mid^2
&\leq&\frac{C}{p_1^2(t)}\|\varphi\|_\infty^2|\bar{u}_0|_H^2.
\end{eqnarray*}
This implies inequality (\ref{3.02}) which completes the proof.
\hfill$\Box$

\end{document}